\documentclass{llncs}
\usepackage{a4wide}
\usepackage{amsfonts}
\usepackage{latexsym}
\usepackage{graphicx}
\usepackage{epsfig}
\usepackage{amsbsy}
\usepackage{amsmath}
\usepackage{amssymb}
\usepackage{psfrag}
\usepackage{subfigure} 
\usepackage{epstopdf}
\usepackage{microtype}

\newcommand{\beq}{\begin{equation}}
\newcommand{\deq}{\end{equation}}
\newcommand{\baq}{\begin{eqnarray}}
\newcommand{\daq}{\end{eqnarray}}

\newcommand{\baqm}{\begin{eqnarray*}}
\newcommand{\daqm}{\end{eqnarray*}}

\def\E{\mathbb{E}}

\usepackage[T1]{fontenc}
\usepackage[utf8]{inputenc}
\usepackage{tabularx,ragged2e,booktabs,caption}
\newcolumntype{C}[1]{>{\Centering}m{#1}}

\begin{document}
\title{Analysis and optimization of vacation and polling models with retrials \protect \footnote{This is an invited, considerably extended version of \cite{BoxmaResing}. The main additions are Subsection \ref{sec:optimization}
and Section \ref{sec:N=N}. These present respectively the optimal behaviour of a single queue system, and the performance analysis for a general
number of queues.}
}
\author{ Murtuza Ali Abidini, Onno Boxma \and Jacques Resing}
\institute{
EURANDOM and Department of Mathematics and Computer Science\\
Eindhoven University of Technology\\
P.O. Box 513, 5600 MB Eindhoven, The Netherlands \\
\email{m.a.abidini@tue.nl,o.j.boxma@tue.nl,j.a.c.resing@tue.nl}}
\maketitle

\begin{abstract}
We study a vacation-type queueing model, and a single-server multi-queue 
polling model, with the special feature of retrials. Just before the server
arrives at a
station there is some deterministic glue period. Customers (both new arrivals and retrials) arriving at
the station during this glue period will be served during the visit of the server. Customers arriving in
any other period leave immediately and will retry after an exponentially distributed time. Our main
focus is on queue length analysis, both at embedded time points (beginnings of glue periods, visit
periods and switch- or vacation periods) and at arbitrary time points.
\\
{\bf Keywords: vacation queue, polling model, retrials}
\end{abstract}

\numberwithin{equation}{section}

\section{Introduction}
Queueing systems with retrials are characterized by the fact that arriving
customers, who find the server busy, do not wait in an ordinary queue.
Instead of that they go into an orbit, retrying to obtain service after
a random amount of time. These systems have received considerable attention
in the literature, see e.g. the book by Falin and Templeton \cite{Falin},
and the more recent book by Artalejo and Gomez-Corral \cite{Artalejo}.

Polling systems are queueing models in which a single server, alternatingly,
visits a number of queues in some prescribed order. Polling systems, too,
have been extensively studied in the literature. For example, various
different service disciplines (rules which describe the server's behaviour
while visiting a queue) and both models with and without switchover times
have been considered. We refer to Takagi \cite{Takagi1,Takagi2} and Vishnevskii and Semenova \cite{Vishnevskii} 
for some literature reviews and to Boon, van der Mei and Winands \cite{Boon}, 
Levy and Sidi \cite{Levy} and Takagi \cite{Takagi3} for overviews of the 
applicability of polling systems.

In this paper, motivated by questions regarding the performance modelling 
of optical networks, we consider vacation and polling systems with retrials.
Despite the enormous amount of literature on both types of models, there are
hardly any papers having both the features of retrials of customers
and of a single server polling a number of queues.
In fact, the authors are only aware of a sequence of papers by 
Langaris \cite{Langaris1,Langaris2,Langaris3} on this topic. In all these papers the author determines
the mean number of retrial customers in the different stations. In \cite{Langaris1} the
author studies a model in which the server, upon polling a station, stays there for an
exponential period of time and if a customer asks for service before this 
time expires, the customer is served and a new exponential stay period at
the station begins. In \cite{Langaris2} the author studies a model with two types of 
customers: primary customers and secondary customers. Primary customers 
are all customers present in the station at the instant the server polls 
the station. Secondary customers are customers who arrive during the 
sojourn time of the server in the station. The server, upon polling a 
station, first serves all the primary customers present and after that 
stays an exponential period of time to wait for and 
serve secondary customers. Finally, in \cite{Langaris3} the author considers a model 
with Markovian routing and stations that could be either of the type of
\cite{Langaris1} or of the type of \cite{Langaris2}.
 
In this paper we consider a polling station with retrials and so-called 
{\em glue} periods. Just before the server arrives at a station there is some
deterministic glue period. Customers (both new arrivals and retrials) 
arriving at the station during this glue period "stick" and will be served 
during the visit of the server. Customers arriving in any other period
leave immediately and will retry after an exponentially distributed time.
 
Our study of queueing systems with retrials
and glue periods was at first instance motivated by questions regarding the performance modelling and analysis
of optical networks.
Optical fibre offers some big advantages for communication w.r.t. copper cables:
huge bandwidth, ultra-low losses, and an extra dimension -- the wavelength of light.
Performance analysis of optical networks is a 
challenging topic (see e.g. Maier \cite{Maier} and Rogiest \cite{Rogiest}). 
In a telecommunication network, packets must be routed from source to destination,
passing through a series of routers and switches. In copper-based transmission links,
packets from different sources are time-multiplexed. This is often modeled
by a single server polling system.
In optical switches, too, one has the need for a protocol to decide which packet may be transmitted.
One might again use a cyclic polling strategy, cyclic meaning that there is a fixed pattern
for giving service to particular ports/stations.
However, unlike electronics, buffering of optical packets is not easy,
as photons cannot be stopped.
Whenever there is a need to buffer photons,
they are made to move locally in fiber loops. These fiber loops or fiber delay lines (FDL) originate and end at the head of a switch.
When a photon arrives at the switch  at a time it cannot be served, it is sent into an FDL, thereby incurring a small delay
to its time of arrival without getting lost or displaced.
Depending on the
availability, requirement, traffic, size of photon and other such factors,
the length (delay produced) of these FDLs can differ. Hence we assume that
these FDLs delay the photons by a random amount of time. Also, if a packet does not receive service after a cycle
through an FDL, then depending on the model
it can go into either the same or a longer or a shorter or randomly to any of
the available FDLs. Hence we assume that two consecutive retrials are
independent of each other.
This FDL feature can be modelled by a retrial queue.

A sophisticated technology that one might try to add to this
is varying the speed of light by changing the refractive index of the fiber loop, cf.\ \cite{okawachi}.
By increasing the refractive index in a small part of
the loop we can achieve `slow light',
which implies slowing the packets.
When a port 
`knows' that it will soon be served,
it may start the process of increasing the refractive index 
at FDLs and at the end of fibers of incoming packets. 
By doing this, it slows down the packets which arrive at the station just before the visit period of the station begins.
This feature is, in our model, incorporated as glue periods immediately before the visit period of the corresponding station.
Packets arriving in this glue period can be served in that subsequent visit period.

The concept of glue period is, to the best of our knowledge, new in polling systems.
It may also be interpreted as a {\em reservation} period.
We view a reservation period as a period in which customers can make a reservation at a station for service in the subsequent visit period
of that station. In our case, such a reservation period occurs immediately before the visit period,
and could be viewed as the last part, $G_i$, of a switchover period of length $S_i+G_i$.
Ordinary gated polling could be viewed as a service discipline in which it is always possible to make a reservation
for the following visit period.

The main contributions of the paper are the following.
(i) For the case of a single queue with vacations and glue periods, we provide a detailed queue-length analysis
at particular embedded epochs and at an arbitrary epoch. We also show how to choose the length of the glue period that
minimizes the mean number of customers in the system.
(ii) We also provide a detailed queue-length analysis for the $N$-queue polling case -- again at particular embedded epochs and
at an arbitrary epoch.

The paper is organized as follows.
In Section~\ref{sec:N=1} we consider the case of a single queue
with vacations and retrials; arrivals and retrials only "stick" during a glue period.
We study this case separately because (i) it is of interest in its own right,
(ii) it allows us to explain the analytic approach
as well as the probabilistic meaning of the main components
in considerable detail, (iii) it makes the analysis
of the multi-queue case more accessible, and (iv) results for the one-queue case may serve as
a first-order approximation for the behaviour of a particular queue in the $N$-queue case,
switchover periods now also representing glue and visit periods at other queues.
In Section~\ref{sec:N=N} the $N$-queue case is analyzed.
Section~\ref{Conclusion} presents some conclusions and suggestions for future research.

\section{Queue length analysis for the single-queue case}
\label{sec:N=1}
\subsection{Model description}
In this section we consider a single queue $Q$ in isolation.
Customers arrive at $Q$ according to a Poisson process with rate $\lambda$.
The service times of successive customers are independent, identically distributed (i.i.d.)
random variables (r.v.), with distribution $B(\cdot)$ and Laplace-Stieltjes transform (LST) $\tilde{B}(\cdot)$.
A generic service time is denoted by $B$.
After a visit period of the server at $Q$ it takes a vacation.
Successive vacation lengths 
are i.i.d. r.v., with $S$ a generic vacation length,
with distribution $S(\cdot)$ and LST $\tilde{S}(\cdot)$.
We make all the usual independence assumptions about interarrival times, service times
and vacation lengths at the queues.
After the server's vacation, a {\em glue} period of deterministic
(i.e., constant) length begins.
Its significance stems from the following assumption.
Customers who arrive at $Q$ do not receive service immediately.
When customers arrive at $Q$ during a glue period $G$, they stick, joining the queue of $Q$.
When they arrive in any other period, they immediately leave and
retry after a retrial interval which is independent
of everything else, and which is exponentially distributed with rate $\nu$.
The glue period is immediately followed by a visit period of the server at $Q$.

The service discipline at $Q$ is gated: During the visit period at $Q$, the server serves all
"glued" customers in that queue, i.e., all customers waiting at the end of the glue period
(but none of those in orbit,
and neither any new arrivals).

We are interested in the steady-state behaviour of this vacation model with retrials.
We hence make the assumption that
$\rho := \lambda \E B < 1$;
it may be verified that this is indeed the condition for this steady-state behaviour to exist.

Some more notation:
\\
$G_{n}$ denotes the $n$th glue period of $Q$.
\\
$V_{n}$ denotes the $n$th visit period of $Q$ (immediately following the $n$th glue period).
\\
$S_{n}$ denotes the $n$th vacation of the server (immediately following the $n$th visit period).
\\
$X_n$ denotes the number of customers 
in the system (hence in orbit) at the start of $G_{n}$.
\\
$Y_{n}$ denotes the number of customers 
in the system at the start of $V_{n}$.
Notice that here we should distinguish between those who are queueing
and those who are in orbit: We write $Y_{n} = Y_{n}^{(q)} + Y_{n}^{(o)}$,
where $q$ denotes queueing and $o$ denotes in orbit.
\\
Finally,
\\
$Z_{n}$ denotes the number of customers
in the system (hence in orbit) at the start of $S_{n}$.

\subsection{Queue length analysis at embedded time points}
\label{queue_lenght_1}
In this subsection we study the steady-state distributions of the numbers of customers at the beginning
of (i) glue periods, (ii) visit periods, and (iii) vacation periods.
Denote by $X$ a r.v. with as distribution the limiting distribution of $X_n$. $Y$ and $Z$ are similarly defined,
and $Y = Y^{(q)} + Y^{(o)}$, the steady-state numbers of customers in queue and in orbit
at the beginning of a visit period (which coincides with the end of a glue period).
In the sequel we shall introduce several generating functions, throughout assuming that their parameter $|z| \leq 1$.
For conciseness of notation, let
$\beta(z) := \tilde{B}(\lambda(1-z))$ and $\sigma(z) := \tilde{S}(\lambda(1-z))$.
Then it is easily seen that
\begin{equation}
\E [z^X] = \sigma(z) \E [z^Z],
\label{eq1}
\end{equation}
since $X$ equals $Z$ plus the new arrivals during the vacation;
\begin{equation}
\E [z^Z] = \E [\beta(z)^{Y^{(q)}} z^{Y^{(o)}}],
\label{eq2}
\end{equation}
since $Z$ equals $Y^{(o)}$ plus the new arrivals during the $Y^{(q)}$ services; and
\begin{equation}
\E [z_q^{Y^{(q)}} z_o^{Y^{(o)}}] = {\rm e}^{-\lambda(1-z_q)G} \E [\{(1 - {\rm e}^{-\nu G})z_q + {\rm e}^{-\nu G} z_o\}^X].
\label{eq3}
\end{equation}
The last equation follows since $Y^{(q)}$ is the sum of new arrivals during $G$ and retrials who return during $G$;
each of the $X$ customers which were in orbit at the beginning of the glue period has a probability $1 - {\rm e}^{-\nu G}$
of returning before the end of that glue period.

Combining Equations (\ref{eq1})-(\ref{eq3}), and introducing
\begin{equation}
f(z) := 
(1 - {\rm e}^{-\nu G})\beta(z) + {\rm e}^{-\nu G} z,
\label{fz}
\end{equation}
we obtain the following functional equation for $\E [z^X]$:
\[
\E [z^X] = \sigma(z) {\rm e}^{-\lambda(1-\beta(z))G} \E [f(z)^X].
\]
Introducing $K(z) := 
\sigma(z) {\rm e}^{-\lambda(1-\beta(z))G}$ and
$X(z) := \E [z^X]$, we have:
\begin{equation}
X(z) = K(z) X(f(z)).
\label{recu}
\end{equation}
This is a functional equation that naturally occurs in the study of queueing models which have a branching-type structure;
see, e.g., \cite{BoxmaCohen} and \cite{Resing93}. Typically, one may view customers who newly arrive into the system during a service as children
of the served customer ("branching"), and customers who newly arrive into the system during a vacation or glue period as immigrants.
Such a functional equation may be solved by iteration, giving rise to an infinite product -- where the $j$th term in the product typically
corresponds to customers who descend from an ancestor of $j$ generations before. In this particular case we have after $n$ iterations:
\begin{equation}
X(z) = \prod_{j=0}^n K(f^{(j)}(z)) X(f^{(n+1)}(z)),
\label{recu1}
\end{equation}
where $f^{(0)}(z) := z$ and $f^{(j)}(z) := f(f^{(j-1)}(z))$, $j=1,2,\dots$.
Below we show that this product converges for $n \rightarrow \infty$ 
iff $\rho < 1$, thus proving the following theorem:
\begin{theorem}
If $\rho < 1$ then the generating function $X(z) = \E [z^X]$ is given by
\begin{equation}
X(z) = \prod_{j=0}^{\infty} K(f^{(j)}(z)).
\end{equation}
\label{thm1}
\end{theorem}
{\bf Proof}.
Equation (\ref{recu}) is an equation for a branching process with 
immigration, where the number of immigrants has generating function $K(z)$
and the number of children in the branching process has generating function
$f(z)$. Clearly, $K^{\prime}(1)= \lambda \E S + \lambda \rho G < \infty$ 
and  $f^{\prime}(1)= {\rm e}^{-\nu G} + \left(1 - {\rm e}^{-\nu G}\right) \rho < 1$, if $\rho < 1$. The result of the theorem now follows directly from the
theory of branching processes with immigration (see e.g., Theorem 1 on 
page $263$ in Athreya and Ney \cite{Athreya}).
\qed
\vspace{0.5cm}
 
\noindent
Having obtained an expression for $\E [z^X]$ in Theorem \ref{thm1},
expressions for
$\E [z^Z]$ and
$\E [z_q^{Y^{(q)}} z_o^{Y^{(o)}}]$
immediately follow from  (\ref{eq2}) and (\ref{eq3}).
Moments of $X$ may be obtained from Theorem \ref{thm1}, but it is also straightforward to obtain
$\E X$ from Equations (\ref{eq1})-(\ref{eq3}):
\begin{eqnarray}
\label{eqn:EZ}
\E X &=& \lambda \E S + \E Z,
\\
\E Z &=& \rho \E Y^{(q)} + \E Y^{(o)},
\\
\E Y^{(q)} &=& \lambda G + (1- {\rm e}^{-\nu G}) \E X,
\\
\E Y^{(o)} &=& {\rm e}^{-\nu G} \E X,
\end{eqnarray}
yielding
\begin{equation}
\E X = \frac{\lambda \E S + \lambda \rho G}{(1-\rho)(1 - {\rm e}^{-\nu G})} .
\label{eqn:EX}
\end{equation}
Hence
\begin{equation}
\E Y^{(q)} = \lambda G + (1- {\rm e}^{-\nu G}) \frac{\lambda \E S + \lambda \rho G}{(1-\rho)(1 - {\rm e}^{-\nu G})} = \frac{\lambda \E S + \lambda G}{1-\rho} ,
\end{equation}
\begin{equation}
\E Y^{(o)} = {\rm e}^{-\nu G} \frac{\lambda \E S + \lambda \rho G}{(1-\rho)(1 - {\rm e}^{-\nu G})} ,
\end{equation}
\begin{equation}
\E Z =  \frac{\lambda \rho G + \lambda \E S [\rho (1- {\rm e}^{-\nu G}) + {\rm e}^{-\nu G}]}{(1-\rho)(1 - {\rm e}^{-\nu G})} .
\end{equation}
Notice that the denominators of the above expressions equal $1-f'(1)$.
Also notice that it makes sense that the denominators contain both the factor $1-\rho$
and the probability $1 - {\rm e}^{-\nu G}$ that a retrial returns during a glue period.

In a similar way as the first moments of $X$, $Y^{(q)}$, $Y^{(o)}$ and $Z$ have been obtained,
we can also obtain their second moment. For further use we here calculate  $\E[ X(X-1)]$:
\begin{eqnarray}
\label{eqn:EXX-1}
\E[ X(X-1)] &=& \frac{K^{\prime\prime}(1)}{(1-\rho)(1 - {\rm e}^{-\nu G})(1 + \rho(1-{\rm e}^{-\nu G}) + {\rm e}^{-\nu G})}
\\
&+&
\frac{K^{\prime}(1)[ 2 K^{\prime}(1) (\rho(1-{\rm e}^{-\nu G}) + {\rm e}^{-\nu G}) + (1-{\rm e}^{-\nu G}) \lambda^2 \E B^2]}{(1-\rho)^2 (1 - {\rm e}^{-\nu G})^2 (1 + \rho(1-{\rm e}^{-\nu G}) + {\rm e}^{-\nu G})},
\nonumber
\end{eqnarray}
where $K^{\prime}(1) = \lambda \E S + \lambda \rho G$ and $K^{\prime\prime}(1) = 
\lambda^2 \E S^2 + 2 \rho \lambda^2 G \E S + \lambda^3 G \E B^2 + (\lambda G \rho)^2$.
\begin{remark}
Special cases of the above analysis are, e.g.:
\\
(i) Vacations of length zero. Simply take $\sigma(z) \equiv 1$ and $\E S=0$ in the above formulas.
\\
(ii) $\nu = \infty$. Retrials now always return during a glue period.
We then have $f(z) = \beta(z)$, which leads to minor simplifications.
\end{remark}
\begin{remark}
It seems difficult to handle the case of non-constant glue periods, as it seems to lead to
a process with complicated dependencies.
If $G$ takes a few distinct values $G_1,\dots,G_N$ with different probabilities, then one might still be able to
obtain a kind of multinomial generalization of the infinite product featuring in Theorem \ref{thm1}.
One would then have several functions $f_i(z) := (1 - {\rm e}^{-\nu G_i}) \beta(z) + {\rm e}^{-\nu G_i} z$,
and all possible combinations of iterations $f_i(f_h(f_k(\dots (z))))$ arising in functions
$K_i(z) := \sigma(z) {\rm e}^{-\lambda(1-\beta(z)) G_i}$, $i=1,2,\dots,N$.
By way of approximation, one might stop the iterations after a certain number of terms, the number depending
on the speed of convergence (hence on $1-\rho$ and on $1 - {\rm e}^{-\nu G_i}$).
\end{remark}
\subsection{Queue length analysis at arbitrary time points}
Having found the generating functions  
of the number of customers at the beginning of (i) glue periods ($\E [z^X]$), (ii) visit
periods $(\E [z_q^{Y^{(q)}} z_o^{Y^{(o)}}])$, and (iii) vacation periods
($\E [z^Z]$), we can also obtain the generating 
function of the number of customers
at arbitrary time points.

\begin{theorem}
If $\rho < 1$, we have the following results: 

\begin{itemize}
\item[a)] The joint generating function, $R_{va}(z_q,z_o)$, of the number of
customers in the queue and in the orbit at an arbitrary time point in a 
vacation period equals the generating function $R_{va}(z_o)$ of the number of customers in orbit
at an arbitrary time point in a vacation period and is given by
\begin{equation}
\label{rva}
R_{va}(z_q,z_o) = \E [z_o^Z]  \frac{1-\tilde{S}(\lambda(1-z_o))}{\lambda(1-z_o) \E S}.
\end{equation}

\item[b)] The joint generating function, $R_{gl}(z_q,z_o)$, of the number of
customers in the queue and in the orbit at an arbitrary time point in a 
glue period is given by
\begin{equation}
R_{gl}(z_q,z_o) = \int_{t=0}^G {\rm e}^{-\lambda(1-z_q)t} 
\E [\{(1 - {\rm e}^{-\nu t})z_q + {\rm e}^{-\nu t} z_o\}^X] \, \frac{dt}{G}.
\end{equation}
\item[c)] The joint generating function, $R_{vi}(z_q,z_o)$, of the number of
customers in the queue and in the orbit at an arbitrary time point in a 
visit period is given by
\begin{equation}
R_{vi}(z_q,z_o) = \frac{z_q\left[\E [z_q^{Y^{(q)}}z_o^{Y^{(o)}}] - 
\E [\tilde{B}(\lambda(1-z_o))^{Y^{(q)}}z_o^{Y^{(o)}}]\right]}{\E [ Y^{(q)}]\left(z_q - \tilde{B}(\lambda(1-z_o))\right)} 
\frac{1-\tilde{B}(\lambda(1-z_o))}{\lambda(1-z_o) \E B} .
\end{equation}
\item[d)] The joint generating function, $R(z_q,z_o)$, of the number of
customers in the queue and in the orbit at an arbitrary time point 
is given by
\begin{equation}
R(z_q,z_o) =  \rho R_{vi}(z_q,z_o) + (1-\rho) \tfrac{G}{G+\E S} 
R_{gl}(z_q,z_o) + (1-\rho) \tfrac{\E S}{G+\E S} R_{va}(z_q,z_o).
\end{equation}
\end{itemize}
\label{arbitrary}
\end{theorem}
{\bf Proof}.
\begin{itemize}
\item[a)] Follows from the fact that during vacation periods there are no customers
in the queue, hence the right-hand side of (\ref{rva}) is independent of $z_q$, and the fact that the number of 
customers at an arbitrary time point in the orbit is the sum of two 
independent terms: The number of customers at the beginning of the 
vacation period and the number that arrived during the past part of the vacation period.
The generating function of the latter is given by
\[
\frac{1-\tilde{S}(\lambda(1-z_o))}{\lambda(1-z_o) \E S}.
\]
\item[b)] Follows from the fact that if the past part of the glue period is equal to 
$t$, the generating function of the number of new arrivals in the queue
during this period is equal to ${\rm e}^{-\lambda(1-z_q)t}$ and  
each customer present in the orbit at the beginning of the glue period
is, independent of the others, still in orbit with probability 
${\rm e}^{-\nu t}$ and has moved to the queue with probability 
$1-{\rm e}^{-\nu t}$.
\item[c)] During an arbitrary point in time in a visit period the number
of customers in the system consists of two parts:
\begin{itemize}
\item the number of customers in the system at the beginning of the service
time of the customer currently in service, leading to the term
(see Remark \ref{pointinvisit} below):
\begin{equation}
\frac{z_q \left(\E [z_q^{Y^{(q)}}z_o^{Y^{(o)}}] - \E [\tilde{B}(\lambda(1-z_o))^{Y^{(q)}}z_o^{Y^{(o)}}]\right)}{\E [Y^{(q)}] \left(z_q - \tilde{B}(\lambda(1-z_o))\right)};
\label{eq:1.3} 
\end{equation}

\item the number of customers that arrived during the past part of the service of the
customer currently in service, leading to the term
\[
\frac{1-\tilde{B}(\lambda(1-z_o))}{\lambda(1-z_o) \E B}.
\]
\end{itemize}
\item[d)] Follows from the fact that the fraction of time the 
server is visiting $Q$ is equal to $\rho$, and if the server is not visiting
$Q$, with probability $\E S/(G+\E S)$ the server is on vacation and with
probability $G/(G+\E S)$ the system is in a glue phase.

\end{itemize}
\qed
\begin{remark}
A straightforward way to prove (\ref{eq:1.3})
is to condition on the number of customers, say, $j$, in queue at the end of a glue period,
and to average the number of customers in queue and in orbit over the beginnings of the $j$ services.
A more elegant proof of Formula \eqref{eq:1.3} uses the observation that 
typically in vacation and polling systems each time a visit beginning or a
service completion occurs, this coincides with a service beginning or a 
visit completion. This observation yields (see, e.g., Boxma, Kella and 
Kosinski \cite{BKK})
\[
\gamma {\cal V}_b(z_q, z_o) + {\cal S}_c(z_q, z_o) = {\cal S}_b(z_q, z_o) 
+\gamma {\cal V}_c(z_q, z_o).
\]
Here, ${\cal V}_b(z_q, z_o)$ and ${\cal V}_c(z_q, z_o)$ are the joint 
generating functions of the number of customers in the queue and in the
orbit at visit beginnings and visit completions, respectively. Similarly,
${\cal S}_b(z_q, z_o)$ and ${\cal S}_c(z_q, z_o)$ are the joint
generating functions of the number of customers in the queue and in the
orbit at service beginnings and service completions, respectively.
Finally, $\gamma$ is the reciprocal of the mean number of customers served 
per visit. 
Clearly,
\[
\gamma = \frac{1}{\E [Y^{(q)}]}, \quad 
{\cal V}_b(z_q, z_o) = \E [z_q^{Y^{(q)}}z_o^{Y^{(o)}}], \quad
{\cal V}_c(z_q, z_o) = \E [\tilde{B}(\lambda(1-z_o))^{Y^{(q)}}z_o^{Y^{(o)}}],
\]
and
\[
{\cal S}_c(z_q, z_o) = \frac{{\cal S}_b(z_q, z_o)}{z_q} \tilde{B}(\lambda(1-z_o)),
\]
which yields that ${\cal S}_b(z_q, z_o)$ is given by \eqref{eq:1.3}.
 \label{pointinvisit}
\end{remark}

From Theorem \ref{arbitrary}, we now can obtain the steady-state mean number of 
customers in the system at arbitrary time points in vacation periods
($\E [R_{va}]$), in glue periods ($\E [R_{gl}]$), in visit periods
($\E [R_{vi}]$) and in arbitrary periods ($\E [R]$). These are given by

\begin{eqnarray}
\E [R_{va}] &=& \E [Z] + \lambda \tfrac{\E [S^2]}{2 \E [S]},\nonumber\\
\E [R_{gl}] &=& \E [X] + \lambda \tfrac{G}{2},\nonumber\\
\E [R_{vi}] &=& 1 + \lambda \tfrac{\E [B^2]}{2 \E [B]} + \tfrac{\E [Y^{(q)} Y^{(o)}]}{\E [Y^{(q)}]} + 
\tfrac{(1+\rho)\E [Y^{(q)}(Y^{(q)}-1)]}{2\E [Y^{(q)}]}, \,\, \nonumber\\
\E [R]      &=& \rho \E [R_{vi}] + (1-\rho) \tfrac{G}{G+\E S} \E [R_{gl}]
+ (1-\rho) \tfrac{\E S}{G+\E S} \E [R_{va}].
\label{eqn:mean}
\end{eqnarray}
Remark that the quantities $\E [Y^{(q)} Y^{(o)}]$ and $\E [Y^{(q)}(Y^{(q)}-1)]$ can be obtained using (\ref{eq3}):
\begin{eqnarray*}
\E [Y^{(q)} Y^{(o)}] &=& \lambda G {\rm e}^{-\nu G} \E [X] + \left(1 -{\rm e}^{-\nu G}\right) {\rm e}^{-\nu G} \E [X(X-1)],\\
\E [Y^{(q)}(Y^{(q)}-1)] &=& (\lambda G)^2 + \left(1 -{\rm e}^{-\nu G}\right)^2 \E [X(X-1)] + 2 \lambda G \left(1 -{\rm e}^{-\nu G}\right) \E [X].
\end{eqnarray*}
By combining these relations with (\ref{eqn:mean}), (\ref{eqn:EZ}), (\ref{eqn:EX}) and (\ref{eqn:EXX-1}),
we obtain -- after tedious calculations
-- the following relatively simple expression for the mean number of customers $\E [R]$:
\begin{equation}
\E [R] = \rho + \frac{\lambda^2 \E [B^2]}{2(1-\rho)}
+ \frac{\lambda \E [(G+S)^2]}{2\E [G+S]} + \frac{\lambda \rho \E [G+S]}{1-\rho}
+ \lambda (\rho G + \E [S]) \frac{{\rm e}^{-\nu G}}{(1-\rho)(1 - {\rm e}^{-\nu G})} ,
\label{EZnew}
\end{equation}
which we rewrite for later purposes as
\begin{eqnarray}
\E [R] &=& \rho + \frac{\lambda^2 \E [B^2]}{2(1-\rho)}
+ \lambda \frac{\E [S]}{G+\E [S]} \frac{\E [S^2]}{2 \E [S]} +
\lambda \frac{\E [S]}{G+\E [S]} G +
\lambda \frac{G}{G+\E [S]} \frac{G}{2}
\nonumber
\\
&+&
\frac{\lambda \rho (G+ \E [S])}{1-\rho}
+ \lambda (\rho G + \E [S]) \frac{{\rm e}^{-\nu G}}{(1-\rho)(1 - {\rm e}^{-\nu G})} .
\label{ER2nd}
\end{eqnarray}

\begin{remark}
(i). It should be noticed that the first two terms in the right-hand side of (\ref{EZnew}) together represent
the mean number of customers in the ordinary $M/G/1$ queue, without vacations and glue periods.
The third term represents the mean number of arrivals during the residual part of a vacation plus glue period.
The fourth term can be interpreted as the mean number of arrivals during a visit period of the server
(since the mean length of one cycle, i.e., visit plus vacation plus glue period, is via a balance argument seen to equal
$\E [C] = \frac{G+\E [S]}{1-\rho}$, while a mean visit period equals $\rho \E [C]$).
The fifth term is the only term involving the retrial rate $\nu$.
In particular, that term disappears when $\nu \rightarrow \infty$. In the latter case, our model reduces to an $M/G/1$ queue with gated vacations,
with vacation lengths $G+S$. The resulting expression for the mean number 
of customers coincides with formula (5.23) of \cite{Takagi4}
(see also formula (3.2.6) of \cite{Tian}). 
\\
(ii). A more interesting limiting operation is to simultaneously let $\nu \rightarrow \infty$ and $G \downarrow 0$,
such that $\nu G$ remains constant.
The resulting model is an $M/G/1$ queue with vacations and binomially gated service; see, e.g., Levy \cite{HLevy}.
Here, each customer who is present at the end of a vacation, will be served in the next visit period
with probability $p = 1- {\rm e}^{-\nu G}$.
In this case, the mean number of customers in the system is given by 
\begin{equation}
\E [R] = \rho + \frac{\lambda^2 \E [B^2]}{2(1-\rho)}
+ \frac{\lambda \E [S^2]}{2\E [S]} + \frac{\lambda \rho \E [S]}{1-\rho}
+ \frac{\lambda \E [S] (1-p)}{p(1-\rho)}.
\end{equation}
This formula coincides with the results obtained in \cite{HLevy} (see 
e.g., formula (7.1) with $N=1$ in \cite{HLevy} for the mean sojourn time 
of customers in this model).
Observe that our formula, after application of Little's formula, does {\it not} 
fully agree with the mean delay expression (5.50b) in \cite{Takagi4}
and with a similar formula on p.\ 90 of
\cite{Tian}. Those mean delay expressions for the binomial gated model 
seem to refer to the case where customers who are not served during a visit (w.p.\ $1-p$) are lost;
hence factors like $\frac{1}{1-p \rho}$ in those mean delay expressions.
\\
(iii). Formula (\ref{ER2nd}) immediately shows how the mean number of customers behaves for very small and for very large values of the glue period length $G$: 
\begin{equation}
\E [R] \sim \frac{\lambda \E [S]}{G \nu (1-\rho)}, ~~~ G \downarrow 0,
\label{G0}
\end{equation}
and
\begin{equation}
\E [R] \sim \frac{\lambda (1+\rho)}{2(1-\rho)} G, ~~~ G \rightarrow \infty .
\label{Ginf}
\end{equation}
We observe that the equations \eqref{G0} and \eqref{Ginf} do not involve $\E[B^2]$
and $\E[S^2]$. Hence the mean number of customers in the system is almost invariant to the variance of
the service times and vacation times when $G$ is either very small or very large.
In Subsection \ref{sec:optimization} we explore the effect of $G$ on $\E [R]$ more deeply.

\end{remark}

\subsection{Optimizing the length of the glue period}
\label{sec:optimization}

One of the main reasons for studying mathematical models of optical 
networks is to improve the performance of the system. In this model the 
length of the glue period is an important system design parameter. The 
results of the previous subsections can, e.g., be used to determine the 
value of $G$ which minimizes the mean number of customers in the system 
at any arbitrary time point. The mean sojourn time of an arbitrary customer
follows from Little's formula. Therefore we can find the value of $G$ 
which minimizes the mean sojourn time of an arbitrary customer.

\begin{figure}
\centering
\subfigure[$\E B=1$]{
\includegraphics[width=0.45\textwidth, height=0.3\textheight]{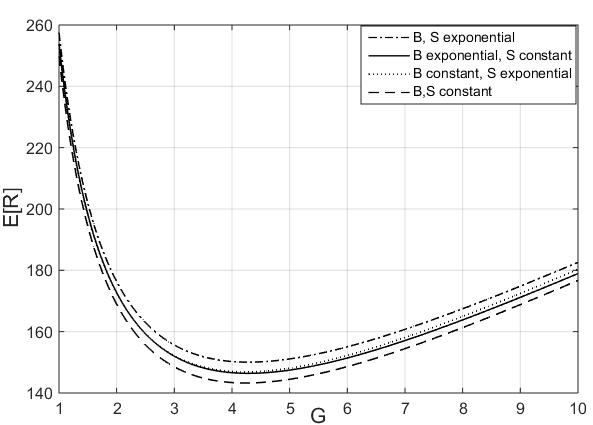} } 
\hspace{5mm}
\subfigure[$\E B=1.1$]{
\includegraphics[width=0.45\textwidth, height=0.3\textheight]{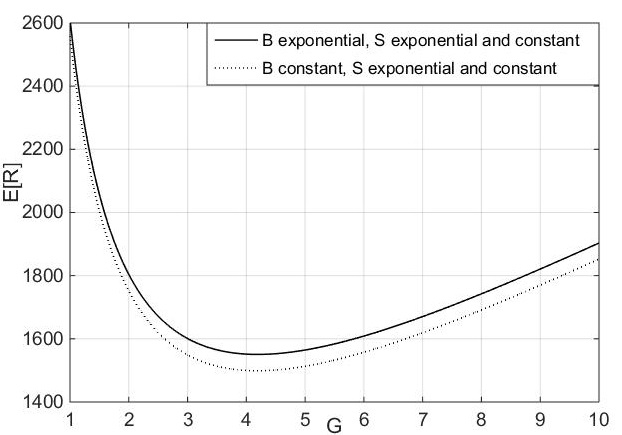}}
\caption{Mean number of customers vs length of glue period}
\label{exp_1}
\end{figure}

Let us first present a sample of numerical results that we obtained for $\E [R]$ as a function of $G$.
We consider four cases: (i) the service time distribution and vacation time distribution are exponential,
(ii) the service time distribution is exponential and the vacation time is constant, (iii) the service time
is constant and the vacation time distribution is exponential and (iv) the service time and vacation
time are constant.
In Fig \ref{exp_1} we take $\lambda = 0.9$, $\nu=0.5$, and $\E S=10$, we plot $G$ vs $\E [R]$, for $\E B =1$
and $1.1.$ Note that $\E B = 1.1$ corresponds to the heavy traffic case $\rho=0.99$. 

From the examples in Fig \ref{exp_1} we observe the following results:
\begin{itemize}
 \item The mean number of customers at any arbitrary point seems to be convex w.r.t. glue period length, i.e., there
 exists a glue period $G_{min}$ where the system has minimum mean number of customers $E[R_{min}]$ and hence
 minimum mean sojourn time.
 \item For a large $G$, $\E[R]$ increases linearly with $G$ (as confirmed by (\ref{Ginf})).
 \item For a very small $G$, $\E[R]$ behaves like $1/G$ (as confirmed by (\ref{G0})).
 \item For different service time and vacation time distributions the mean number of customer
 is changed but the difference in $G_{min}$ is insignificantly small.
 \item As $\rho \rightarrow 1$ , i.e. the system is in heavy traffic regime, the value of $\E[R]$
 becomes insensitive to the distribution of $S$.
 \end{itemize}

Indeed, if $G$ is very small, the number of customers joining the queue in each glue period is very small
and thereby the efficiency of the system is decreased. On the other hand, a large $G$ means the system stays in the glue period for
a long time and this decreases the efficiency of the system. Hence it is logical to have a $G_{min}$ which optimizes the system.
 
We will now prove that $E[R]$ is indeed convex in $G$.
Twice differentiating the expression for $\E [R]$ in Equation \eqref{EZnew} w.r.t.\ $G$ gives
 \begin{equation}
\frac{d^2}{dG^2}\E [R] = \frac{\lambda}{\E [G+S]} \Bigg[\frac{ \E [(G+S)^2]}{\E [G+S]^2} -1 \Bigg] + \frac{\lambda \nu {\rm e}^{-\nu G} }{(1-\rho)(1- {\rm e}^{-\nu G})^2}
\Bigg[\nu(\rho G + \E [S]) \frac{ (1 + {\rm e}^{-\nu G})}{(1 - {\rm e}^{-\nu G})}  -2 \rho \Bigg]. \\
\label{DG2}
\end{equation}
Clearly, the first term in the right-hand side of (\ref{DG2}) is nonnegative. 
Let $$Q(G) := \nu(\rho G + \E [S]) \frac{ (1 + {\rm e}^{-\nu G})}{(1 - {\rm e}^{-\nu G})}  -2 \rho .$$
We can see $ Q(G) \longrightarrow \infty$ as $G \longrightarrow 0$ or $G \longrightarrow \infty$.
Let $Q(G)$ attain its minimum at $G = g.$
Hence calculating the first derivative of $Q(G)$ and equating it to zero, at $G=g,$ gives
$$ \rho g + \E [S] = \frac{\rho}{2 \nu {\rm e}^{-\nu g}} (1-{\rm e}^{-2 \nu g}) .$$
Therefore $$Q(g) =  \frac{\rho}{2} [{\rm e}^{\nu g} -2 + {\rm e}^{-\nu g}] \geq 0.$$
We observe that the minimum value of $Q(G)$ is always nonnegative.
Since both terms of (\ref{DG2}) are nonnegative, $\frac{d^2}{dG^2}\E [R] \geq 0.$

Hence $\E [R]$, the mean number of customers at an arbitrary point of time in the system, is convex in $G$.
So the system can improve the quality of service by setting an optimal value of $G$ for the fixed glue period. 

In Table \ref{tab:1} we analyze the behaviour of $G_{min}$ and $\E[R_{min}]$ as we increase $\E [S]$
for an exponential distribution $B(\cdot)$ with $\E B=1$, arrival rate $\lambda=0.5$ 
and retrial rate $\nu=0.5$.
\begin{table}[h]
\resizebox{\textwidth}{!}{
\begin{minipage}{\linewidth}
\begin{center}
\captionof{table}{Optimal glue period length and corresponding minimum mean number of customers}
\vspace{0.5cm}
\label{tab:1} 
\begin{tabular}{| l | c  |c | c |r|}
     \hline
    $\E S$ & $G_{min} $ for & $G_{min}$ for  & $\E[R_{min}]$ for &$\E[R_{min}]$ for  \\
     &   exponential $S(\cdot)$ &  constant $S(\cdot)$  &   exponential $S(\cdot)$  & constant $S(\cdot)$ \\ \hline
    $0$ &  $\epsilon$, ~~$\epsilon \longrightarrow 0$  &$\epsilon$, ~~$\epsilon \longrightarrow 0$  & $ \sim 2$ &   $\sim 2$\\ \hline
    0.1 & 0.608& 0.606 & 2.476  & 2.473 \\ \hline    
    0.5 & 1.334 & 1.320& 3.419   & 3.385  \\ \hline
    1 & 1.846 & 1.801& 4.259  & 4.170\\ \hline
    5 & 3.694 & 3.508& 9.276  & 8.549 \\ \hline
    10 & 4.785& 4.495 & 14.778  & 13.070 \\ \hline
    50 & 7.783 & 7.225& 56.027  & 45.156  \\ \hline
    100 & 9.167& 8.521 & 106.606  & 83.634 \\ \hline
  
      \end{tabular} 
\end{center}
      
\end{minipage}}
\end{table}

Table \ref{tab:1} suggests that, in the case under consideration,
\begin{itemize}
 \item $G_{min}$ and $\E[R_{min}]$ increase when $\E [S]$ increases.
 \item $G_{min}$ and $\E[R_{min}]$ increase when the variance of $S$ becomes larger.
 \item When $\E S$ approaches $0$, $G_{min}$ also approaches $0$.
When there is no customer in the queue, the system will then have a series of very short
glue periods, and when a customer arrives or returns from orbit, it can almost instantaneously be taken into service.
In this case, the system reduces to an ordinary $M/G/1$ retrial queue;
indeed, Formula (\ref{EZnew}) reduces to $\E [R] = \rho + \frac{\lambda^2 \E [B^2]}{2(1-\rho)}$ $ + \frac{\lambda \rho}{\nu (1-\rho)}$
which is in agreement with Formula (1.37) of \cite{Falin}.
\end{itemize}

\section{Queue length analysis for the N-queue case}
\label{sec:N=N}
\subsection{Model description}
In this section we consider a one-server polling model
with multiple queues,
$Q_i,$ $i= 1,\cdots,N$.
Customers arrive at $Q_i$ according to a Poisson process with rate $\lambda_i$;
they are called type-$i$ customers,
$i=1,\cdots,N$.
The service times at $Q_i$ are i.i.d. r.v., with $B_i$ denoting a generic service time,
with distribution $B_i(\cdot)$ and LST $\tilde{B}_i(\cdot)$,
$i=1,\cdots,N$.
The server follows cyclic polling, thus after a visit of $Q_i$, it switches to $Q_{i+1}$.
Successive switchover times from $Q_i$ to $Q_{i+1}$
are i.i.d. r.v., with $S_i$ a generic switchover time,
with distribution $S_i(\cdot)$ and LST $\tilde{S}_i(\cdot)$,
$i=1,\cdots,N$.
We make all the usual independence assumptions about interarrival times, service times
and switchover times at the queues.
After a switch of the server to $Q_i$, there first is a deterministic (i.e., constant)
glue period $G_i$, before the visit of the server at $Q_i$ begins, $i=1,\cdots,N$.
As in the one-queue case, the significance of the glue period stems from the following assumption.
Customers who arrive at $Q_i$ do not receive service immediately.
When customers arrive at $Q_i$ during a glue period $G_i$, they stick, joining the queue of $Q_i$.
When they arrive in any other period, they immediately leave and
retry after a retrial interval which is independent
of everything else, and which is exponentially distributed with rate $\nu_i$, $i=1,\cdots,N$.

The service discipline at all queues is gated: During the visit period at $Q_i$, the server serves all
"glued" customers in that queue, i.e., all type-$i$ customers waiting at the end of the glue period --
but none of those in orbit,
and neither any new arrivals.

We are interested in the steady-state behaviour of this polling model with retrials.
We hence assume that the stability condition 
$\sum_{i=1}^N \rho_i < 1$ holds, where
$\rho_i := \lambda_i \E B_i$.

Some more notation:
\\
$G_{ni}$ denotes the $n$th glue period of $Q_i$.
\\
$V_{ni}$ denotes the $n$th visit period of $Q_i$.
\\
$S_{ni}$ denotes the $n$th switch period out of $Q_i$, $i=1,\cdots,N$.
\\
$(X_{n1}^{(i)},X_{n2}^{(i)},\cdots ,X_{nN}^{(i)})$ denotes the vector of numbers of customers of type $1$ to type $N$
in the system (hence in orbit) at the start of $G_{ni}$, $i=1,\cdots,N$.
\\
$(Y_{n1}^{(i)},Y_{n2}^{(i)},\cdots ,Y_{nN}^{(i)})$ denotes the vector of numbers of customers of type $1$ to type $N$
in the system at the start of $V_{ni}$, $i=1,\cdots,N$.
We distinguish between those who are queueing in $Q_i$
and those who are in orbit for $Q_i$: We write $Y_{ni}^{(i)} = Y_{ni}^{(iq)} + Y_{ni}^{(io)}$, 
$i=1,\cdots ,N$, where $q$ denotes queueing and $o$ denotes in orbit.
\\
Finally,
\\
$(Z_{n1}^{(i)},Z_{n2}^{(i)},\cdots ,Z_{nN}^{(i)})$ denotes the vector of numbers of customers of type $1$ to type $N$
in the system (hence in orbit) at the start of $S_{ni}$, $i=1,\cdots ,N$.

\subsection{Queue length analysis}
In this subsection we study the steady-state joint distribution of the
numbers of customers in the system at beginnings of glue periods.
This will also immediately yield the steady-state joint distributions of the numbers of customers
in the system at the beginnings of visit periods and of switch periods.
We follow a similar generating function approach as in the one-queue case, throughout making the following assumption regarding
the parameters of the generating functions: $|z_i| \leq 1$, $|z_{iq}| \leq 1$,
$|z_{io}| \leq 1$.
Observe that the generating function of the vector of numbers of arrivals
at $Q_1$ to $Q_N$ during a type-$i$ service time $B_i$ is $\beta_i(z_1,z_2,\cdots ,z_N) := \tilde{B}_i(\sum_{j=1}^{N}\lambda_j(1-z_j))$.
Similarly, 
the generating function of the vector of numbers of arrivals
at $Q_1$ to $Q_N$ during a type-$i$ switchover time $S_i$ is $\sigma_i(z_1,z_2,\cdots ,z_N) := \tilde{S}_i(\sum_{j=1}^{N}\lambda_j(1-z_j))$.
Since the server serves $Q_{i+1}$ after $Q_{i}$ we can successively express
(in terms of generating functions) $(X_{n1}^{(i+1)},X_{n2}^{(i+1)},\cdots ,X_{nN}^{(i+1)})$ into $(Z_{n1}^{(i)},Z_{n2}^{(i)},\cdots ,Z_{nN}^{(i)})$,
$(Z_{n1}^{(i)},Z_{n2}^{(i)},\cdots ,Z_{nN}^{(i)})$ into $(Y_{n1}^{(i)},Y_{n2}^{(i)},\cdots ,Y_{ni}^{(iq)}, Y_{ni}^{(io)},\cdots ,Y_{nN}^{(i)})$,
and \\ $(Y_{n1}^{(i)},Y_{n2}^{(i)},\cdots ,Y_{ni}^{(iq)}, Y_{ni}^{(io)},\cdots ,Y_{nN}^{(i)})$ into $(X_{n1}^{(i)},X_{n2}^{(i)},\cdots ,X_{nN}^{(i)})$; etc.
Denote by $(X_{1}^{(i)},X_{2}^{(i)},\cdots ,X_{N}^{(i)})$ the vector with as distribution the limiting distribution of
$(X_{n1}^{(i)},X_{n2}^{(i)},\cdots ,X_{nN}^{(i)})$, $i=1,\cdots ,N$, and similarly introduce
$(Z_1^{(i)},Z_2^{(i)},\cdots ,Z_N^{(i)})$ and
$(Y_1^{(i)},Y_2^{(i)},\cdots ,Y_N^{(i)})$,
with $Y_i^{(i)} = Y_i^{(iq)} + Y_i^{(io)}$ ,
for $i=1,\cdots ,N$.
We have: 
\begin{equation}
\E \left[z_1^{X_{1}^{(i+1)}} z_2^{X_{2}^{(i+1)}}\cdots z_N^{X_{N}^{(i+1)}}\right] =
\sigma_i(z_1,z_2,\cdots,z_N) 
\E \left[z_1^{Z_{1}^{(i)}} z_2^{Z_{2}^{(i)}}\cdots z_N^{Z_{N}^{(i)}}\right] .
\label{XZ}
\end{equation}
\begin{eqnarray*}
\hspace{-0.2in}
\E \left[z_{1}^{Z_{1}^{(i)}} z_2^{Z_{2}^{(i)}}\cdots z_N^{Z_{N}^{(i)}}| Y_{1}^{(i)} = h_{1}, Y_{2}^{(i)} = h_{2},\cdots,Y_{i}^{(iq)} = h_{iq}, Y_{i}^{(io)} = h_{io},\cdots, Y_{N}^{(i)} = h_N\right]\nonumber\\
=\left(\prod_{j\ne i}z_j^{h_j}\right) z_i^{h_{io}} \left[\beta_i(z_1,z_2,\cdots,z_N)\right]^{h_{iq}},
\end{eqnarray*}
yielding
\begin{equation}
\E \left[z_{1}^{Z_{1}^{(i)}} z_2^{Z_{2}^{(i)}}\cdots z_N^{Z_{N}^{(i)}}\right]
=
\E \left[ \left[\beta_i(z_1,z_2,\cdots,z_N)\right]^{Y_{i}^{(iq)}} z_1^{Y_{1}^{(i)}}z_2^{Y_{2}^{(i)}}\cdots z_i^{Y_{i}^{(io)}}\cdots z_N^{Y_{N}^{(i)}}\right] .
\label{ZY}
\end{equation}
Furthermore,
\begin{eqnarray*}
\E \left[ z_1^{Y_{1}^{(i)}}z_2^{Y_{2}^{(i)}}\cdots z_{iq}^{Y_{i}^{(iq)}}z_{io}^{Y_{i}^{(io)}}\cdots z_N^{Y_{N}^{(i)}}
| X_{1}^{(i)} = a_1, X_{2}^{(i)} = a_2,\cdots,X_{N}^{(i)} = a_N\right]
\nonumber
\\
=
\left(\prod_{j\ne i}
z_j^{a_j} {\rm e}^{-\lambda_j(1-z_j) G_i}\right) {\rm e}^{-\lambda_i(1-z_{iq}) G_i}
\left[(1 - {\rm e}^{- \nu_i G_i}) z_{iq} + {\rm e}^{-\nu_i G_i} z_{io}\right]^{a_i} ,
\end{eqnarray*}
yielding
\begin{eqnarray}
\hspace{-0.3in}
\E \left[z_1^{Y_{1}^{(i)}}z_2^{Y_{2}^{(i)}}\cdots z_{iq}^{Y_{i}^{(iq)}} 
z_{io}^{Y_{i}^{(io)}}\cdots z_N^{Y_{N}^{(i)}}\right]
&=&
\left(\prod_{j\ne i} {\rm e}^{-\lambda_j(1-z_j) G_i}\right) {\rm e}^{-\lambda_i(1-z_{iq}) G_i} \nonumber \\
 &\times& \E \left[ \left(\prod_{ j\ne i}z_j^{X_j^{(i)}}\right)\left[(1 - {\rm e}^{- \nu_i G_i}) z_{iq} + {\rm e}^{-\nu_i G_i} z_{io}\right]^{X_{i}^{(i)}}\right].
\label{YX}
\end{eqnarray}
It follows from (\ref{XZ}), (\ref{ZY}) and (\ref{YX}), with
\begin{equation*}
f_i(z_1,z_2,\cdots ,z_N) := (1-{\rm e}^{-\nu_i G_i}) \beta_i(z_1,z_2,\cdots ,z_N) + {\rm e}^{-\nu_i G_i} z_i,
\end{equation*}
that
\begin{eqnarray}
\hspace*{-0.3in}
\E \left[z_1^{X_{1}^{(i+1)}} z_2^{X_{2}^{(i+1)}}\cdots z_N^{X_{N}^{(i+1)}}\right] &=&
\sigma_i(z_1,z_2,\cdots ,z_N)\left(\prod_{ j\ne i} {\rm e}^{-\lambda_j(1-z_j) G_i}\right)
{\rm e}^{-\lambda_i(1-\beta_i(z_1,z_2,\cdots ,z_N)) G_i}
\nonumber \\
&\times& \E \left[ \left(\prod_{ j\ne i}z_j^{X_j^{(i)}}\right)\left[f_{i}(z_1,z_2,\cdots ,z_N)\right]^{X_{i}^{(i)}}\right] .
\label{een}
\end{eqnarray}
Let $z=(z_1,z_2,\cdots,z_N)$;
further
\begin{eqnarray*}
 h_{i}(z)&:=&f_i(z_1,\cdots ,z_i,h_{i+1}(z),\cdots ,h_N(z)), \\
 h_N(z)&:=&f_N(z_1,\cdots ,z_N),\\
 \beta^{(i)}(z) &:=& \beta_i(z_1,\cdots ,z_i,h_{i+1}(z),\cdots ,h_N(z)), \\
  \beta^{(N)}(z) &:=& \beta_N (z_1,\cdots ,z_N),\\ 
  \sigma^{(i)}(z) &:=& \sigma_i(z_1,\cdots ,z_i,h_{i+1}(z),\cdots ,h_N(z)), \\
  \sigma^{(N)}(z) &:=& \sigma_N (z_1,\cdots ,z_N).
 \end{eqnarray*}
Since the server moves to $Q_1$ after $Q_N$, substituting $i=N$ in \eqref{een}, we have 
\begin{eqnarray}
\E \left[z_1^{X_{1}^{(1)}} z_2^{X_{2}^{(1)}}\cdots z_N^{X_{N}^{(1)}}\right]& =&
\sigma^{(N)}(z)\left(\prod_{ j\ne N} {\rm e}^{-\lambda_j(1-z_j) G_N}\right) {\rm e}^{-\lambda_N\left(1-\beta^{(N)}(z)\right) G_N} \nonumber \\
&\times&  \E \left[ \left(\prod_{ j\ne N}z_j^{X_j^{(N)}}\right)\left[h_N(z)\right]^{X_{N}^{(N)}}\right] .
\label{X1N} 
\end{eqnarray}
From \eqref{een} we have 
\begin{eqnarray}
\hspace*{-0.3in}
 \E \left[ \left(\prod_{ j\ne N}z_j^{X_j^{(N)}}\right)[h_N(z)]^{X_{N}^{(N)}} \right]&=& 
\sigma^{(N-1)}(z)\left(\prod_{ j=1}^{N-2} {\rm e}^{-\lambda_j(1-z_j) G_{N-1}}\right) \nonumber \\
&\times& {\rm e}^{-\lambda_{N-1}(1-\beta^{(N-1)}(z)) G_{N-1}} {\rm e}^{-\lambda_N(1-h_N(z)) G_{N-1}} \nonumber \\
&\times& \E \left[\left(\prod_{ j=1}^{N-2}z_j^{X_j^{(N-1)}}\right)[h_{N-1}(z)]^{X_{N-1}^{(N-1)}}[h_N(z)]^{X_{N}^{(N-1)}} \right]. 
\label{XNN-1}
\end{eqnarray}

From \eqref{X1N} and \eqref{XNN-1} we have
\begin{eqnarray*}
\hspace*{-0.2in}
 \E \Big[z_1^{X_{1}^{(1)}} && z_2^{X_{2}^{(1)}}\cdots z_N^{X_{N}^{(1)}}\Big] =
\sigma^{(N)}(z)\left(\prod_{ j\ne N} {\rm e}^{-\lambda_j(1-z_j) G_N}\right) {\rm e}^{-\lambda_N(1-\beta^{(N)}(z)) G_N} \nonumber \\
 \times && \sigma^{(N-1)}(z)\left(\prod_{ j=1}^{N-2} {\rm e}^{-\lambda_j(1-z_j) G_{N-1}}\right) {\rm e}^{-\lambda_{N-1}(1-\beta^{(N-1)}(z)) G_{N-1}} \nonumber \\
 \times && {\rm e}^{-\lambda_N(1-h_N(z)) G_{N-1}}\E \left[\left(\prod_{ j=1}^{N-2}z_j^{X_j^{(N-1)}}\right)[h_{N-1}(z)]^{X_{N-1}^{(N-1)}}[h_N(z)]^{X_{N}^{(N-1)}} \right]. \nonumber \\
=&&\left(\prod_{i=N-1}^{N}\sigma^{(i)}(z) {\rm e}^{-G_i D_i(z)}\right)\E \left[\left(\prod_{ j=1}^{N-2}z_j^{X_j^{(N-1)}}\right)[h_{N-1}(z)]^{X_{N-1}^{(N-1)}}[h_N(z)]^{X_{N}^{(N-1)}}  \right], \nonumber \\
 \end{eqnarray*}
 where
 \begin{equation*}
  D_i(z)= \sum_{j=1}^{i-1} \lambda_j (1-z_j) +\lambda_i \left(1-\beta^{(i)}(z)\right) + \sum_{j=i+1}^{N} \lambda_j(1-h_j(z)).
 \end{equation*}
By recursively substituting as above we get
\begin{equation}
 \E \left[z_1^{X_{1}^{(1)}} z_2^{X_{2}^{(1)}}\cdots z_N^{X_{N}^{(1)}}\right]=\prod_{i=1}^{N}\sigma^{(i)}(z) \prod_{i=1}^{N} {\rm e}^{-G_i D_i(z)} 
 \E \left[ [h_1(z)]^{X_{1}^{(1)}}[h_2(z)]^{X_{2}^{(1)}}\cdots[h_N(z)]^{X_{N}^{(1)}}\right].
 \label{M1}
 \end{equation}

Equation \eqref{M1} can be divided into three factors, representing the switchover period, glue period and
visit period respectively. The first factor, for a particular value of $i$, represents the arrivals during
the switchover time after the visit of $Q_i$. The second factor represents the arrivals during the glue period
before a visit of $Q_i$. It is further divided into three generating functions.
First are the arrivals of type $j <i$; these don't have any further effect on the system. Then the arrivals of type $i$, these are served during the following
visit and produce new children (i.e., arrivals during their service) of each type. Finally those of type $j>i$ which may or may not be served in future visits and if served produce new children of each type.
These two factors are taken for all $i=1,\cdots,N$. The third factor represents the descendants (arrivals during services, arrivals
during services of customers who arrived during services, etc.) of $(X_1^{(1)},\cdots,X_N^{(1)})$.
 
Consider 
\begin{equation*}
 X(z)=\E \left[\left(\prod_{j=1}^{N}z_j^{X_{j}^{(1)}}\right)\right],
\end{equation*}
with an obvious definition of $K(z)$, we can rewrite \eqref{M1} into
\begin{equation}
 X(z)=K(z)X(h(z))
 \label{Branching}
\end{equation}
where 
\begin{equation*}
 h(z):=(h_{1}(z),\cdots ,h_N(z)).
\end{equation*}
Define, for all $i=1,\cdots, N$,
\[
h_i^{(0)}(z)= z_i, \quad h_i^{(n)}(z)= h_i(h_1^{(n-1)}(z),h_2^{(n-1)}(z),\cdots,h_N^{(n-1)}(z)).
\]

\begin{theorem}
 If $\sum_i \rho_i < 1$, then the generating function  $ X(z)$ is given by 
 \begin{equation}
 X(z)=\prod_{m=0}^{\infty}K(h_1^{(m)}(z),h_2^{(m)}(z),\cdots,h_N^{(m)}(z)) .
 \label{chainbranch}
 \end{equation}
 \label{thm:chainbranch}
\end{theorem}

{\bf Proof.} 
Equation \eqref{chainbranch} follows from \eqref{Branching} by iteration.
We still need to prove that the infinite product converges if $\sum_i \rho_i < 1$.
Equation (\ref{Branching}) is an equation for a multi-type branching process with 
immigration, where the number of immigrants of different types has 
generating function $K(z)$ and the number of children of different 
types of a type $i$ individual in the branching process has generating 
function $h_i(z)$, $i=1,\cdots,N$. An important role in the analysis of such a process is 
played by the mean matrix $M$ of the branching process,
\begin{equation}
M = \left(\begin{array}{cc}
m_{11} \cdots m_{1N} \\
\vdots ~~~\ddots ~~~\vdots\\
m_{N1} \cdots m_{NN}
\end{array} \right),
\label{meanmatrix}
\end{equation}
where $m_{ij}$ represents the mean number of children of type $j$ of a type $i$ individual.
The elements of the matrix $M$ are the same as given in Section 5 of Resing \cite{Resing93}, which is 
\begin{equation}
m_{ij}= f_{ij} \cdot 1[j \leq i] + \sum_{k=i+1}^{N} f_{ik} m_{kj},
\end{equation}
where
$m_{ij}= \frac{\partial h_i}{\partial z_j}(1,1,\cdots,1)$
 and
 $f_{ij}= \frac{\partial f_i}{\partial z_j}(1,1,\cdots,1).$

 We observe that the equation for $m_{ij}$ is the sum of two terms. First the children of type $j \leq i$,
 who do not affect the system in the future. Next the children of type $j$ produced by the children of type $k>i$ in the subsequent visits.

The theory of multi-type branching processes with immigration (see Quine 
\cite{Quine} and Resing \cite{Resing93}) now states that if (i) the expected 
total number of immigrants in a generation is finite and (ii) the maximal 
eigenvalue $\lambda_{max}$ of the mean matrix $M$ satisfies 
$\lambda_{max} < 1$, then the generating function of the steady state 
distribution of the process is given by (\ref{chainbranch}). 
To complete the proof of Theorem \ref{thm:chainbranch}, we shall now verify (i) and (ii).
\\

\noindent
Ad (i):
The expected total number of 
immigrants in a generation is 
\begin{eqnarray}
\lambda_1 \Bigg(G_1 \Bigg(\sum_j m_{1j}\Bigg) + \sum_{j=1}^N \E S_j + \sum_{j=2}^N G_j \Bigg)~~~~~~&& \nonumber \\
+ \lambda_2 \Bigg(\big((G_1+ \E S_1)\big(1- {\rm e}^{-\nu_2G_2}\big) + G_2\big) \Bigg(\sum_j m_{2j}&\Bigg)& + (G_1+\E S_1){\rm e}^{-\nu_2G_2}+ \sum_{j=2}^N \E S_j + \sum_{j=3}^N G_j \Bigg) \nonumber \\
+ \cdots ~~~~~~~~~~~~~~~~~~~~~~~~~~~~~~~~~~~~~~~~~~~~~~~~~~~~~~ & &\nonumber \\
+ \lambda_N \Bigg(\Bigg(\Bigg(\sum_{j=1}^{N-1}(G_j+\E S_j)\Bigg)\big(1- {\rm e}^{-\nu_N G_N}\big) + G_N&\Bigg)& \Bigg(\sum_j m_{Nj}\Bigg) + \sum_{j=1}^{N-1}(G_j+\E S_j){\rm e}^{-\nu_N G_N}+  \E S_N \Bigg) \nonumber \\
= \sum_i\lambda_i \Bigg(\Bigg(\Bigg(\sum_{j=1}^{i-1}(G_j+\E S_j)\Bigg)\big(1- {\rm e}^{-\nu_i G_i}\big) + G_i&\Bigg)& \Bigg(\sum_j m_{ij}\Bigg)+\sum_{j=1}^{i-1}(G_j+\E S_j){\rm e}^{-\nu_i G_i}+ \sum_{j=i}^N \E S_j + \sum_{j=i+1}^N G_j \Bigg). \nonumber \\
\label{immigrants}
\end{eqnarray}
Since the above equation is a finite sum/product of finite terms it is indeed finite.

Here, the term $\lambda_1 (G_1 (\sum_j m_{1j}))$ corresponds to
the type $1$ customers arriving during the glue period of $Q_1$ and their subsequent children of all types.
The term $\lambda_1(\sum_{j=1}^N \E S_j + \sum_{j=2}^N G_j )$ corresponds to the type $1$ customers
arriving during the glue periods of $Q_j$, $j=2,\cdots,N,$ and switchover periods after
$Q_j$, $j=1,\cdots,N$. These customers arrive after the visit of $Q_1$ and hence do not get served or 
produce children.
The term $\lambda_2 (((G_1+\E S_1)(1- {\rm e}^{-\nu_2G_2}) + G_2) (\sum_j m_{2j})+ (G_1+\E S_1){\rm e}^{-\nu_2G_2})$
corresponds to the type $2$ customers arriving during the glue period of $Q_1,~Q_2$,
the switchover period after $Q_1$ and their subsequent children. The term
$\lambda_2(\sum_{j=2}^N \E S_j + \sum_{j=3}^N G_j )$ corresponds to the type $2$ customers arriving during the
glue periods of $Q_j$, $j=3,\cdots,N,$ and switchover periods after $Q_j$, $j=2,\cdots,N$. These customers
do not produce any children. Similarly  the term
$\lambda_N (((\sum_{j=1}^{N-1}(G_j+\E S_j))(1- {\rm e}^{-\nu_N G_N}) + G_N) (\sum_j m_{Nj}) + \sum_{j=1}^{N-1}(G_j+\E S_j){\rm e}^{-\nu_N G_N})$
corresponds to the type $N$ customers arriving during the glue period of $Q_1, \cdots, Q_N$,
the switchover periods after $Q_1,\cdots,Q_{N-1}$ and their subsequent children. The term
$\lambda_N \E S_N$ corresponds to the type $N$ customers arriving during the
switchover period after $Q_N$, which do not produce any children.
\\
\noindent
Ad (ii): Define the matrix

\begin{equation}
H=\begin{pmatrix}
{\rm e}^{-\nu_1G_1} + \left(1-{\rm e}^{-\nu_1G_1}\right)\rho_1 &~& \left(1-{\rm e}^{-\nu_1G_1}\right)\lambda_2 \E B_1&~& \cdots &~& \left(1-{\rm e}^{-\nu_1G_1}\right)\lambda_N \E B_1\\\
\left(1-{\rm e}^{-\nu_2G_2}\right)\lambda_1 \E B_2 &~&{\rm e}^{-\nu_2G_2} + \left(1-{\rm e}^{-\nu_2G_2}\right) \rho_2&~& \cdots &~&\left(1-{\rm e}^{-\nu_2G_2}\right)\lambda_N \E B_2 \\\
\vdots &~&\vdots &~&\cdots &~&\vdots\\
\left(1-{\rm e}^{-\nu_NG_N}\right)\lambda_1 \E B_N &~& \left(1-{\rm e}^{-\nu_NG_N}\right)\lambda_2\E B_N &~&\cdots &~&{\rm e}^{-\nu_NG_N} + \left(1-{\rm e}^{-\nu_NG_N}\right) \rho_N
\end{pmatrix},
\label{meanmatrix}
\end{equation}
where the elements $h_{ij}$ of the matrix $H$ represent the mean number of type $j$ customers that replace a type $i$ 
customer during a visit period of $Q_i$ (either new arrivals if the customer is served, or the
customer itself if it is not served).
We have that  
\begin{equation}
H 
\begin{pmatrix}
\E B_1 \\
\E B_2 \\
\vdots \\
\E B_N
\end{pmatrix}
= \begin{pmatrix}
\left[{\rm e}^{-\nu_1G_1} + \left(1-{\rm e}^{-\nu_1G_1}\right) (\sum_j\rho_j )\right] \E B_1
\\
\left[{\rm e}^{-\nu_2G_2} + \left(1-{\rm e}^{-\nu_2G_2}\right) (\sum_j \rho_j)\right] \E B_2
\\
\vdots \\
\left[{\rm e}^{-\nu_NG_N} + \left(1-{\rm e}^{-\nu_NG_N}\right) (\sum_j \rho_j)\right] \E B_N
\end{pmatrix} 
<
\begin{pmatrix}
\E B_1 \\
\E B_2 \\
\vdots \\
\E B_N
\end{pmatrix}\end{equation}
if and only if $\sum_j\rho_j < 1$. Using this result and 
following the same line of proof as in 
Section 5 of Resing \cite{Resing93}, we can show that the stability 
condition $\sum_j\rho_j < 1$ implies that also the 
maximal eigenvalue $\lambda_{max}$ of the mean matrix $M$ satisfies 
$\lambda_{max} < 1$. This concludes the proof.
\qed

\vspace{0.5in}

We can now obtain the moments, $\E X_{j}^{(i+1)}$, either from \eqref{chainbranch} or in a similar way as in Section \ref{queue_lenght_1},
in terms of $\E X_{j}^{(i)}$ and $\E X_{i}^{(i)}$ :

\begin{equation*}
\E X_{j}^{(i+1)} = \lambda_{j} \E S_{i} + \E Z_{j}^{(i)}.
\end{equation*}
When $j \ne i$,
\begin{equation*}
\E Z_{j}^{(i)} = \lambda_{j} \E B_{i} \E Y_{i}^{(iq)} + \E Y_{j}^{(i)}, 
\end{equation*}
else
\begin{equation*}
\E Z_{i}^{(i)} = \lambda_{i} \E B_{i} \E Y_{i}^{(iq)} + \E Y_{i}^{(io)}.
\end{equation*}
Further
\begin{eqnarray*}
\E Y_j^{(i)} &=& \lambda_j G_i +  \E X_j^{(i)},
\\
\E Y_{i}^{(iq)} &=& \lambda_i G_i + (1- {\rm e}^{-\nu_i G_i}) \E X_i^{(i)},
\\
\E Y_i^{(io)} &=& {\rm e}^{-\nu_i G_i} \E X_i^{(i)}.
\end{eqnarray*}
From the above equations we get, when $j \ne i$ :
\begin{equation*}
\E X_{j}^{(i+1)} = \lambda_{j} \E S_{i} + \lambda_{j} (1+\rho_i) G_i+ \lambda_{j} \E B_{i} (1- {\rm e}^{-\nu_i G_i}) \E X_i^{(i)} + \E X_{j}^{(i)},
\end{equation*}
and
\begin{equation*}
\E X_{i}^{(i+1)} = \lambda_{i} \E S_{i} + \lambda_{i} \rho_i G_i+ (\rho_i(1- {\rm e}^{-\nu_i G_i}) + {\rm e}^{-\nu_i G_i})\E X_i^{(i)} .
\end{equation*}

Using flow balance arguments (mean number of customers of type $i$ served per cycle equals mean number of type $i$ customers
arriving per cycle) and the obvious fact that the mean cycle time equals $\E C := \sum_i (\E S_i + G_i) / (1-\rho)$ , we obtain
\begin{equation}
\E Y_{i}^{(iq)} =   \frac{\lambda_i}{1- \rho}\sum_j (\E S_j + G_j).
\label{eqn:y_work}
\end{equation}
We can also use a similar argument for mean number of type $i$ customers leaving the orbit, $(1-{\rm e}^{-\nu_i G_i}) \E X_i^{(i)}$,
to equal the mean number of type $i$ customers entering it, $\lambda_i (\E C - G_i)$, per cycle, yielding
\begin{equation}
 \E X_{i}^{(i)} =  \frac{\lambda_i}{1-{\rm e}^{-\nu_i G_i}}[\frac{\sum_j (\E S_j + G_j)}{1- \rho} - G_i].
\label{eqn:x_work}
 \end{equation}
 We can observe that \eqref{eqn:y_work} and \eqref{eqn:x_work} satisfy the above relation between $\E Y_{i}^{(iq)}$ and $\E X_{i}^{(i)}$.
Further for each $i$, cyclically substituting we get all $\E X_{j}^{(i)}$ and therefore $\E Y_{j}^{(i)}$ and $\E Z_{j}^{(i)}$.

The second moments of $X_{j}^{(1)}$ and the various terms $\E [X_{j}^{(i)} X_{k}^{(i)}]$ can be obtained
by solving a set of equations which is derived by twice differentiating \eqref{Branching} w.r.t. $z_j$ and $z_k$, $j,k = 1,\cdots,N$, and calculating
the value at $z =(1,1,\cdots,1)$. Since the system is cyclic, once we obtain $\E {X_j^{(1)}}^2$, $j=1,\cdots,N$, we can
similarly obtain $\E {X_j^{(i)}}^2$, $i	=1,\cdots,N$, by changing indices.
It is not difficult to develop an efficient procedure for determining higher moments
in polling systems with a branching discipline, cf.\ \cite{Resing93}.

\subsection{Queue length analysis at arbitrary time points}

In the previous section we have given the procedure for finding the distribution of the number of customers at the beginning
of (i) glue periods ($\E [\prod_j z_j^{X_j^{(i)}}]$), (ii) visit
periods $(\E [\left(\prod_{j\ne i} z_j^{Y_j^{(i)}}\right)z_{iq}^{Y^{(iq)}} z_{io}^{Y^{(io)}}])$, and (iii) switchover periods
($\E [\prod_j z_j^{Z_j^{(i)}}]$), for $i=1,\cdots,N$. Similar to the vacation model, we now obtain the generating 
functions of the numbers of customers in queue and in the orbit, at all the stations, at arbitrary time points.

\begin{theorem}
If $\sum_i\rho_i < 1$ and $(z_q,z_o) :=(z_{1q},z_{1o},\cdots,z_{Nq},z_{No})$,  we have the following results: 

\begin{itemize}
\item[a)] The joint generating function, $R^{(i)}_{sw}(z_{q},z_{o})$, of the numbers of
customers in the queue and in the orbit at an arbitrary time point in a 
switchover period after $Q_i$ equals the joint generating function, $R^{(i)}_{sw}(z_o)$,
of the numbers of customers in orbit at an arbitrary time point in a 
switchover period after $Q_i$ and is given by
\begin{equation}
R^{(i)}_{sw}(z_q,z_o) = \E [\prod_jz_{jo}^{Z_j^{(i)}}]  \frac{1-\tilde{S_i}(\sum_j\lambda_j(1-z_{jo}))} {\left(\sum_j\lambda_j(1-z_{jo})\right) \E S_i}.
\end{equation}
\item[b)] The joint generating function, $R^{(i)}_{gl}(z_q,z_o)$, of the numbers of
customers in the queue and in the orbit at an arbitrary time point in a 
glue period of $Q_i$ is given by
\begin{equation}
\hspace{-0.3in}
R^{(i)}_{gl}(z_q,z_o) = \left(\prod_{j \ne i} \E [z_{jo}^{X_j^{(i)}}]\right) \int_{t=0}^{G_i} \left(\prod_{j\ne i}{\rm e}^{-\lambda_j(1-z_{jo})t}\right){\rm e}^{-\lambda_i(1-z_{iq})t} 
\E [\{(1 - {\rm e}^{-\nu_i t})z_{iq} + {\rm e}^{-\nu_i t} z_{io} \} ^{X_i^{(i)}}] \, \frac{dt}{G_i}.
\end{equation}
\item[c)] The joint generating function, $R^{(i)}_{vi}(z_q,z_o)$, of the numbers of
customers in the queue and in the orbit at an arbitrary time point in a 
visit period of $Q_i$ is given by
\begin{eqnarray}
\hspace{-0.3in}
R^{(i)}_{vi}(z_q,z_o) &=& \frac{z_{iq} \left(\E [z_{iq}^{Y_i^{(iq)}}\left(\prod_{j\ne i}z_{jo}^{Y_j^{(i)}} \right)z_{io}^{Y_i^{(io)}}] - \E [\tilde{B_i}(\sum_j \lambda_j(1-z_{jo}))^{Y_i^{(iq)}}\left(\prod_{j\ne i}z_{jo}^{Y_j^{(i)}} \right)z_{io}^{Y_i^{(io)}}]\right)}{\E [Y_i^{(iq)}] \left(z_{iq} - \tilde{B_i}(\sum_j \lambda_j(1-z_{jo}))\right)} \nonumber \\ 
&& \times \frac{1-\tilde{B_i}(\sum_j\lambda_j(1-z_{jo}))} {\left(\sum_j\lambda_j(1-z_{jo})\right) \E B_i}.
\end{eqnarray}
\item[d)] The joint generating function, $R(z_q,z_o)$, of the numbers of
customers in the queue and in the orbit at an arbitrary time point 
is given by
\begin{equation}
\hspace{-0.3in}
R(z_q,z_o) = \sum_i\left( \rho_i R^{(i)}_{vi}(z_q,z_o) + (1-\rho_i) \tfrac{G_i}{\sum_j(G_j+\E S_j)} 
R^{(i)}_{gl}(z_q,z_o) + (1-\rho_i) \tfrac{\E S_i}{\sum_j(G_j+\E S_j)} R^{(i)}_{sw}(z_q,z_o)\right).
\end{equation}

\end{itemize}
\label{arbitraryN}
\end{theorem}
{\bf Proof}.
The proof follows the same lines as the proof of Theorem \ref{arbitrary}, in particular for parts a and d.
We restrict ourselves here to an outline of the proof of parts b and c.
\begin{itemize}
\item[b)] Follows from the fact that if the past part of the glue period is equal to 
$t$, the generating function of the number of new arrivals  of type $i$ in the queue
during this period is equal to ${\rm e}^{-\lambda_i(1-z_{iq})t}$ and  
each  type $i$ customer present in the orbit at the beginning of the glue period
is, independent of the others, still in orbit with probability 
${\rm e}^{-\nu_i t}$ and has moved to the queue with probability 
$1-{\rm e}^{-\nu_i t}$. Further the generating function of the number
of new arrivals  of any type $j \ne i$ in the queue
during this period is equal to ${\rm e}^{-\lambda_j(1-z_{jo})t}.$
\item[c)] During an arbitrary point in time in a visit period the number
of customers in the system consists of two parts:
\begin{itemize}
\item the number of customers in the system at the beginning of the service
time of the customer currently in service, leading to the term
\begin{equation*}
\hspace{-0.5in}
\frac{z_{iq} \left(\E [z_{iq}^{Y_i^{(iq)}}\left(\prod_{j\ne i}z_{jo}^{Y_j^{(i)}} \right)z_{io}^{Y_i^{(io)}}] - \E [\tilde{B_i}(\sum_j \lambda_j(1-z_{jo}))^{Y_i^{(iq)}}\left(\prod_{j\ne i}z_{jo}^{Y_j^{(i)}} \right)z_{io}^{Y_i^{(io)}}]\right)}{\E [Y_i^{(iq)}] \left(z_{iq} - \tilde{B_i}(\sum_j \lambda_j(1-z_{jo}))\right)};
\end{equation*}
 
(see Remark \ref{pointinvisit}). 
\item the number of customers that arrived during the past part of the service of the
customer currently in service, leading to the term
\[
\frac{1-\tilde{B_i}(\sum_j\lambda_j(1-z_{jo}))} {\left(\sum_j\lambda_j(1-z_{jo})\right) \E B_i}.
\]
\end{itemize}
\qed
\end{itemize}

From Theorem \ref{arbitraryN}, we now can obtain the steady-state mean number of 
customers in the system at arbitrary time points in switchover periods
($\E [R^{(i)}_{sw}]$) after $Q_i$, in glue periods ($\E [R^{(i)}_{gl}]$) and in visit periods
($\E [R^{(i)}_{vi}]$) of $Q_i$ , for $i=1,\cdots,N$, and at any arbitrary time point ($\E [R]$). These are given by

\begin{eqnarray}
\E [R^{(i)}_{sw}] &=& \sum_j(\E [Z_j^{(i)}] + \lambda_j \tfrac{\E [S_i^2]}{2 \E [S_i]}),\nonumber\\
\E [R^{(i)}_{gl}] &=& \sum_j(\E [X_j^{(i)}] + \lambda_j \tfrac{G_i}{2}),\nonumber\\
\E [R^{(i)}_{vi}] &=& 1 + (\sum_j\lambda_j) \tfrac{\E [B_i^2]}{2 \E [B_i]} + \tfrac{\E [Y_i^{(iq)} Y_i^{(io)}]}{\E [Y_i^{(iq)}]} + (\sum_{j \ne i} \tfrac{\E [Y_i^{(iq)} Y_j^{(i)}]}{\E [Y_i^{(iq)}]}) +
\tfrac{(1+\E B_i\sum_j \lambda_j)\E [Y_i^{(iq)}(Y_i^{(iq)}-1)]}{2\E [Y_i^{(iq)}]}, \,\, \nonumber\\
\E [R]      &=& \sum_i \left(\rho_i \E [R^{(i)}_{vi}] + (1-\rho) \tfrac{G_i}{\sum_j(G_j+\E S_j)} \E [R^{(i)}_{gl}]
+ (1-\rho) \tfrac{\E S_i}{\sum_j(G_j+\E S_j)} \E [R^{(i)}_{sw}]\right).
\label{meannum_N}
\end{eqnarray}
The mean number of type $k$ customers in the system at arbitrary time points in a switchover period after $Q_i$ and a glue
period before $Q_i$ are given by the values of the $k$-th term of the sums in the formulas of $\E [R^{(i)}_{sw}]$ and $\E [R^{(i)}_{gl}]$.
The mean number of type $i$ customers in the system at arbitrary time points in a visit period of $Q_i$ is given by
$$1 + \lambda_i \tfrac{\E [B_i^2]}{2 \E [B_i]} + \tfrac{\E [Y_i^{(iq)} Y_i^{(io)}]}{\E [Y_i^{(iq)}]} +
\tfrac{(1+\rho_i)\E [Y_i^{(iq)}(Y_i^{(iq)}-1)]}{2\E [Y_i^{(iq)}]} .$$

The quantities $\E [Y_i^{(iq)} Y_i^{(io)}]$, $\E [Y_i^{(iq)} Y_j^{(i)}]$ and $\E [Y_i^{(iq)}(Y_i^{(iq)}-1)]$ can be obtained using (\ref{YX}).


\begin{remark}
Using a similar approach as presented in \cite{HLevy} for a
polling system with binomial-gated service, we can also obtain the 
following expression for $\E [R_i]$, the steady-state mean number of type-$i$ 
customers in the system at arbitrary time points,
\begin{equation}
\E [R_i ] = \rho_i +  
\tfrac{\E [Y_i^{(iq)} Y_i^{(io)}]}{\E [Y_i^{(iq)}]} + 
\tfrac{(1+\rho_i)\E [Y_i^{(iq)}(Y_i^{(iq)}-1)]}{2\E [Y_i^{(iq)}]}, 
\label{meannum_levy}
\end{equation}
which, after summing over $i$, leads to the alternative formula
\begin{equation}
\E [R ] = \rho + \sum_{i=1}^N \tfrac{\E [Y_i^{(iq)} Y_i^{(io)}]}
{\E [Y_i^{(iq)}]} +
\sum_{i=1}^N \tfrac{(1+\rho_i)\E [Y_i^{(iq)}(Y_i^{(iq)}-1)]}{2\E [ Y_i^{(iq)}]}. 
\label{mean_levy}
\end{equation}
\end{remark}

\begin{remark}
From (\ref{mean_levy}), we can derive an explicit expression for the mean 
number of customers in the system in the case of 
a completely symmetric system ($\lambda_i = \lambda/N$, $\rho_i = \rho/N$, $\E B_i = \E B$,
$\E B_i^2 = \E B^2$, $\E S_i = \E S$, $\E S_i^2 = \E S^2$, $G_i =G$, $\nu_i = \nu$). In this case 
we get

\begin{eqnarray}
\E [R] &=& \rho + \frac{\lambda^2 \E [B^2]}{2(1-\rho)}
+ \frac{\lambda N (G + \E S)}{2} + \frac{\lambda {\rm Var}(S)}{2(G+ \E S)} \nonumber \\
&+& \frac{(N+1) \lambda \rho (G + \E S)}{2(1-\rho)}
+ \frac{\lambda {\rm e}^{-\nu G}}{1 - {\rm e}^{-\nu G}}\left[\frac{N(G + \E S)}{1-\rho} - G\right].
\label{meannum_N1}
\end{eqnarray}
\end{remark}
\begin{remark}
\label{rm7}
In \cite{Boxma89} the following so-called pseudo conservation law  -- an explicit expression
for $\sum \rho_i \E[W_i]$, with $\E[W_i]$ the mean waiting time of a customer
of type $i$ until the start of its service -- has been proven for a large class
of polling systems, which also contains the present model:
  \begin{equation}
   \sum_i \rho_i \E [W_i] = \rho \frac{\sum_i \lambda_i \E [B_i^2]}{2(1-\rho)} + \rho \frac{\E[(\sum_i (S_i + G_i))^2]}{2\E [\sum_i (S_i + G_i)]} + \frac{\E [\sum_i (S_i + G_i)]}{2(1-\rho)}\left[\rho^2-\sum_i \rho_i^2 \right] + \sum_i \E [F_i],
   \label{pcl_1}
  \end{equation}
where $\sum_i (S_i + G_i)$ is the sum  of all the idle periods of the server and $F_i$ is the work left in $Q_i$ at the start of a switchover from $Q_i$.
Hence, $ \E[F_i] = \E[Z_i^{(i)}] \E[B_i]$. Other than $\E[F_i]$, the expression is independent of the service discipline.
Using a fairly straightforward mean value analysis we obtain 
 \begin{equation}
 \E[F_i] =\rho_i^2 \E [C] + \frac{\rho_i}{{\rm e}^{\nu_i G_i}-1} \left( \E [C] - G_i\right).
\label{moment1}
\end{equation}

 From Equations \eqref{pcl_1} and \eqref{moment1}, we obtain the following pseudo conservation law:
 \begin{eqnarray}
 \sum_i \rho_i \E W_i &=& \rho \left[\frac{\sum_i \lambda_i \E [B_i^2]}{2(1-\rho)}+ \frac{\E [(\sum_i (S_i +  G_i))^2]}{2 \E [\sum_i (S_i +  G_i)]}\right] 
+ \frac{\E [\sum_i (S_i +  G_i)]}{2(1-\rho)}\left[\rho^2 + \sum_i \rho_i^2 \right] \nonumber \\
&&+ \sum_i \frac{\rho_i}{{\rm e}^{\nu_i G_i}-1}\left[\frac{\E [\sum_i (S_i +  G_i)]}{1-\rho} - G_i \right].
\label{pcl_2}
\end{eqnarray}
The use of this pseudo conservation law seems to be the easiest
way to derive (\ref{meannum_N1}).
Another useful aspect of this pseudo conservation law is that it allows us to study the effect of
the length of the glue period. 
\\
In a more general (not necessarily optical switch related) setting, referred to in the Introduction,
the glue period may represent the only opportunity to make a reservation for service.
The glue or reservation period now is the last part of the switchover period to $Q_i$; one could view $S_i^* := S_i + G_i$
as the total switchover period into $Q_i$ and $G_i$ as the reservation period.
Formula (\ref{pcl_2}) allows us to study how the mean waiting times or mean queue lengths
are affected by having only a brief reservation period, instead of being able to make a reservation at any time
(which is the classical gated polling system).
The only difference
in $\sum \rho_i \E[W_i]$
between our reservation model and the classical gated
polling model
is the last term.
If $\nu_i G_i$ is very small, that last term will be dominant.
If it is not very small whereas, e.g., the second moments of the service times are large, then
the extra term is relatively small -- and the advantage for the system operator of having to offer only very limited reservation opportunities may outweigh
the fact that waiting times and queue lengths become slightly larger.
\end{remark}

\vspace{0.5cm}
\subsection{Numerical example}
\label{numerical}

In this subsection we present some numerical results for the $2$-queue case.
We have numerically evaluated the expressions for $\E[R_1]$, $\E[R_2]$ and $\E[R]$ using Eq. \eqref{meannum_N}.
We have also verified that Eq. \eqref{meannum_levy} gives exactly the same $\E[R_i]$ values,
and that the pseudo conservation law \eqref{pcl_2} 
is satisfied in the numerical examples.
In the numerical example in this section we consider a model with two 
stations. We have chosen $\lambda_1=2$, $\lambda_2=1$, $\nu_1=1$ and $\nu_2=0.2$.
The switchover times are deterministic with $\E[S_1]=0.5$ and $\E[S_2]=1$.
The service times are exponential with $\E[B_1]=\E[B_2]=0.2$.

\begin{figure}
\centering
\includegraphics[width=0.45\textwidth, height=0.25\textheight]{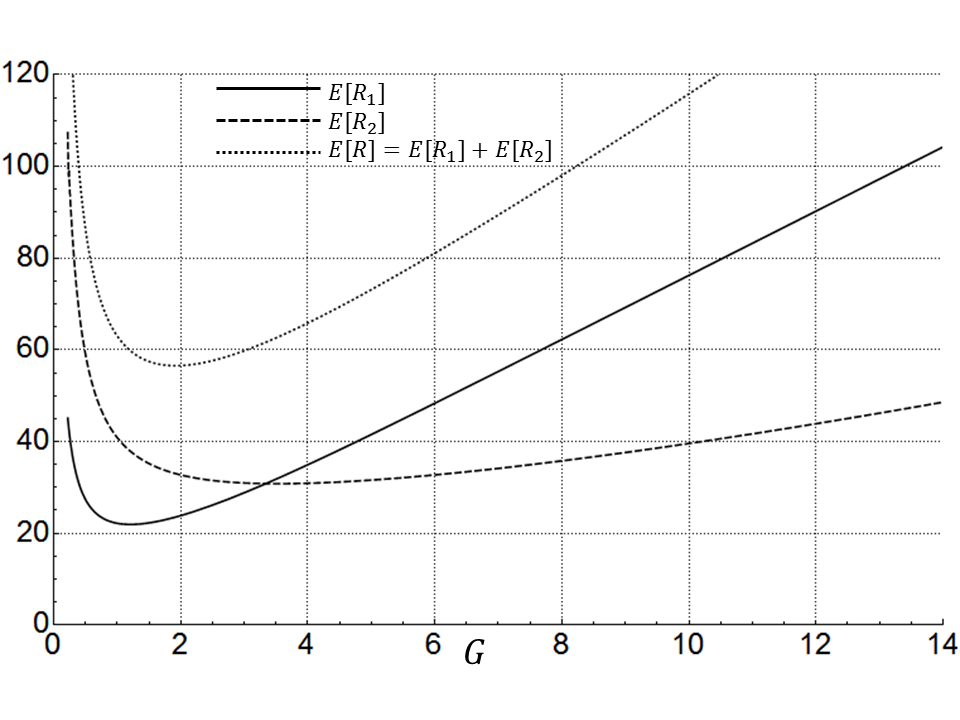} 
\caption{Mean number of customers vs length of glue period in case $G_1=G_2=G$.}
\label{equalG}
\end{figure}

First of all we look at the mean number of type-$1$ customers in the system, $\E[R_1]$, the mean number of type-$2$ customers in the system, $\E[R_2]$,
and the mean total number of customers in the system, $\E[R]=\E[R_1]+\E[R_2]$, 
if we vary the length of the glue period $G$ in the case that $G_1=G_2=G$
(see Fig \ref{equalG}).

\begin{figure}[h]
\centering
\subfigure[varying $G_1$, with $G_2=2$]{
\includegraphics[width=0.45\textwidth, height=0.25\textheight]{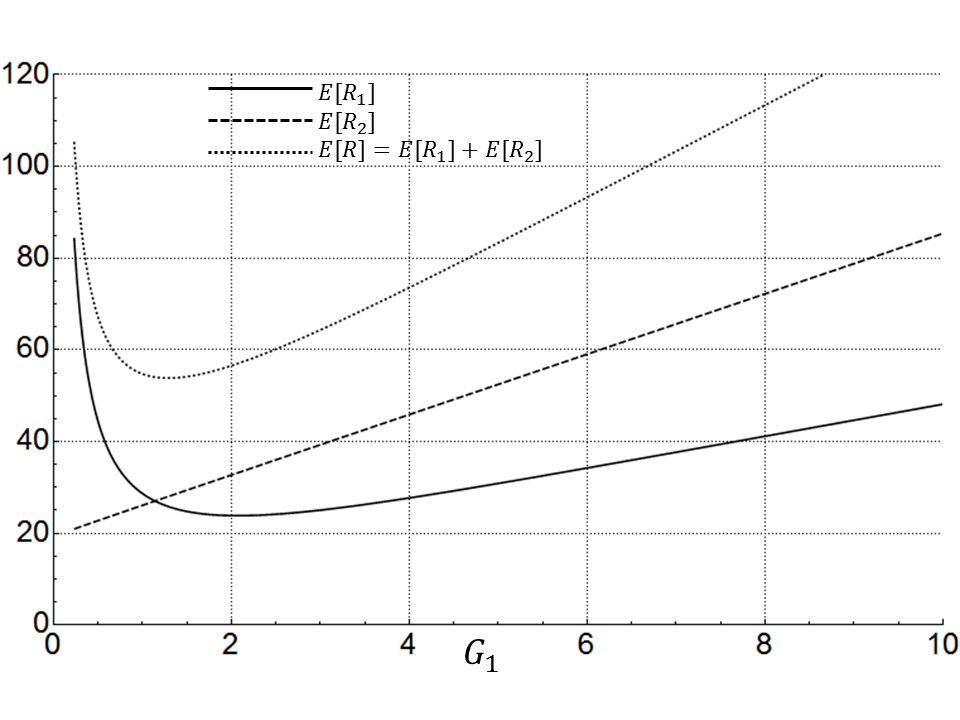}} 
\hspace{5mm}
\subfigure[varying $G_2$, with $G_1=2$]{
\includegraphics[width=0.45\textwidth, height=0.25\textheight]{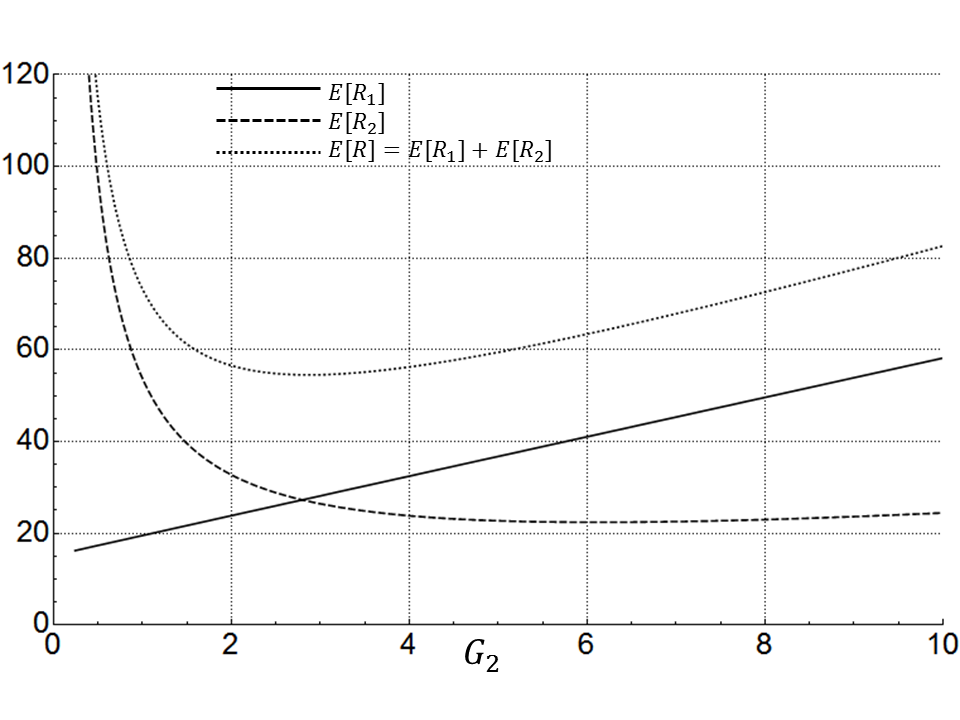}}
\caption{Mean number of customers vs length of glue period in case $G_1 \neq G_2$.}
\label{unequalG}
\end{figure}

Next we look in Fig \ref{unequalG} at the case where the glue periods 
$G_1$ and $G_2$ are unequal. In particular we look at the mean number of 
customers in the system of type-$1$, type-$2$ and in total if we vary one 
of the glue periods while keeping the other glue period constant. 

The plots suggest that there is again a unique value of the length of the 
glue period that minimizes the mean total number of customers in the 
system. They also show that the length of the glue period that is optimal 
for $\E[R]$ is different from the one that is optimal for $\E[R_1]$ and 
the one that is optimal for $\E[R_2]$.

\vspace{0.5cm}
\section{Conclusions and suggestions for future research}
\label{Conclusion}

In this paper we have studied vacation queues and $N$-queue polling models
with the gated service discipline and
with retrials.
Motivated by optical communications,
we have introduced a glue period just before a server visit;
during such a glue period,
new customers and retrials "stick" instead of immediately going into orbit.
For both the vacation queue and the $N$-queue polling model,
we have derived steady-state queue length distributions at an arbitrary epoch and at various specific epochs.
This was accomplished by establishing a relation to branching processes.

Jointly with B. Kim and J. Kim, in \cite{KimJac} we are considering the model with non-constant glue periods.
 We derive $\E [R]$ using MVA (mean value analysis) for the vacation model. It gives
 us the result for general glue period distributions without calculating the generating
 functions at different epochs. We show that a workload decomposition and
pseudo conservation law, as discussed in \cite{Boxma89} and Remark ~\ref{rm7}, can be derived for these variants
and generalizations, and they may be exploited for analysis and optimization purposes.
We shall then also try to explore the following observation:
One may view our $N$-queue model as a polling model with a new variant of binomial gated,
with adaptive probability $p_i$ of serving a customer at a visit of $Q_i$; $p_i=1$ when the customer arrived
in the preceding glue period, and $p_i = 1 - {\rm e}^{-\nu_i G_i}$ otherwise.
We would also like to explore the possibility to study the heavy traffic behavior 
of these models via the relation to branching processes, cf.\ \cite{RDM}.

From a more applied perspective we are looking into systems with multiple customer classes where
one class has priority over another. This would help to incorporate the real life scenario
where some type of data packets, like video buffering,
should have as low latency as possible, whereas others, like a file transfer,
can be delayed a bit longer.

Finally, we would like to point out an important advantage of optical fibre:
the wavelength of light. A fibre-based network node may thus route incoming packets not only
by switching in the time-domain, but also by wavelength division multiplexing.
In queueing terms, this gives rise to {\em multiserver} polling models,
each server representing a wavelength. We refer to
\cite{Antunes1} for the stability analysis of multiserver polling models,
and to \cite{Antunes2} for a mean field approximation of large passive optical networks.
Therefore we would like to study multiserver polling models
with the additional features of retrials and glue periods.
\\

\noindent
{\bf Acknowledgment}
\\

The authors gratefully acknowledge fruitful discussions with
Kevin ten Braak and Tuan Phung-Duc
about retrial queues and with Ton Koonen about optical networks.
The research is supported by the IAP program BESTCOM, funded by the Belgian government,
and by the Gravity program NETWORKS, funded by the Dutch government.

\end{document}